\documentclass[12pt]{article}
\usepackage{theorem,amsfonts}

\newtheorem{theorem}{Theorem}[section]
\newtheorem{proposition}{Proposition}[section]
\newtheorem{lemma}{Lemma}[section]
\newtheorem{corollary}{Corollary}[section]
\theorembodyfont{\upshape}

\newtheorem{remark}{Remark}
\newtheorem{proof}{Proof}

\newtheorem{acknowledgement}{Acknowledgement}

\newcommand{\bt}{\begin{theorem}}
\newcommand{\et}{\end{theorem}}
\newcommand{\bl}{\begin{lemma}}
\newcommand{\el}{\end{lemma}}
\newcommand{\bp}{\begin{proposition}}
\newcommand{\ep}{\end{proposition}}
\newcommand{\bo}{\begin{proof}}
\newcommand{\eo}{\end{proof}}
\newcommand{\br}{\begin{remark}}
\newcommand{\er}{\end{remark}}
\newcommand{\bc}{\begin{corollary}}
\newcommand{\ec}{\end{corollary}}
\newcommand{\be}{\begin{enumerate}}
\newcommand{\ee}{\end{enumerate}}

\title{Krengel-Lin decomposition for probability measures on hypergroups}
\author{C. Robinson Edward Raja}
\date{}

\begin{document}
\maketitle

\let\epsi=\epsilon
\let\vepsi=\varepsilon
\let\lam=\lambda
\let\Lam=\Lambda 
\let\ap=\alpha
\let\vp=\varphi
\let\ra=\rightarrow
\let\Ra=\Rightarrow 
\let\LRa=\Leftrightarrow
\let\Llra=\Longleftrightarrow
\let\Lla=\Longleftarrow
\let\lra=\longrightarrow
\let\Lra=\Longrightarrow
\let\ba=\beta
\let\ga=\gamma
\let\Ga=\Gamma
\let\un=\upsilon

\begin{abstract}
A Markov operator $P$ on a $\sigma$-finite measure space $(X, \Sigma , m)$
with invariant measure $m$ is said to have Krengel-Lin decomposition 
if $L^2 (X) = E_0 \oplus L^2 (X,\Sigma _d)$ where 
$E_0 = \{ f \in L^2 (X) \mid ||P^n (f) || \ra 0 \}$ 
and $\Sigma _d$ is the deterministic $\sigma $-field of $P$.  
We consider convolution operators and we show
that a measure $\lam$ on a hypergroup has Krengel-Lin decomposition if and
only if the sequence $(\check \lam ^n *\lam ^n )$ converges to an
idempotent or $\lam$ is scattered.  We verify this condition for
probabilities on Tortrat groups, on commutative hypergroups and on central
hypergroups.  We give a counter-example to show that the decomposition is
not true for measures on discrete hypergroups.
\end{abstract}


\begin{section}{Introduction and generalities}

Let $K$ be a locally compact space.  Let ${\cal P}(K)$ be the space of 
regular Borel probability measures on $K$ with the weak topology that is the
smallest topology for which the functions $f \mapsto \mu (f)$ from
${\cal P} (K)$ to $\mathbb R$ are continuous for any bounded continuous 
function on $K$.  For any $x \in K$, let $\delta _x$ be the probability
measure concentrated at the point $x$ and for any $\mu \in {\cal P}(K)$,
${\rm supp} (\mu )$ is the support of $\mu$ (which is smallest closed
subset
$C$ of $K$ for which $\mu (C) =1$).  

A locally compact space $K$ with
a binary operation $*$ on the space $M^b(K)$ of bounded measures
on $K$ is called a {\it hypergroup} if 

\be

\item $*$ is bilinear and separately continuous from $M^b(K)
\times M^b(K) \ra M^b(K)$ so that $(M^b(K), +, *)$ is an 
associative algebra where $+$ is the usual additive operation,  

\item the mapping $(\mu , \lam ) \mapsto \mu *\lam$ from 
${\cal P}(K) \times {\cal P}(K)$ into ${\cal P}(K)$ is continuous, 

\item for every $x, y \in K$, the support of $\delta _x *\delta _y$ is
compact and the mapping $(x, y) \mapsto {\rm supp} 
(\delta _x*\delta _y )$
from $K\times K$ to ${\cal A}(K)$, the space of
compact subsets of $K$ with the Michael topology is continuous,

\item there exists an element $e \in K$ such that 
$\delta _e *\delta _x = \delta
_x *\delta _e = \delta _x$ for all $x \in K$, 

\item there exists an involution $x \mapsto \check x$ of $K$
such that $e$ is contained in the support of $\delta _x *\delta _y$ if and
only if $x = \check y$ and $\mu \mapsto \check \mu$ is an
anti-homomorphism of ${\cal P}(K)$. 
\ee 

We assume that 
$K$ is a $\sigma$-compact hypergroup with a right-invariant Haar
measure $m$.  It is known that locally compact groups, 
commutative hypergroups, discrete hypergroups and central hypergroups
admit invariant measures.  For any locally compact group $G$ and any
compact subgroup $M$ of $G$, the double coset space $G//M$ with
the quotient topology and the convolution induced by the group operation
in $G$ is a hypergroup with an invariant measure (which is induced by the
Haar measure on $G$): see \cite {B-H} for results on hypergroups.  

Let $L^2(K)$ be the space of square integrable function on $K$ with 
respect to $m$.  For $f \in L^2 (K)$ and for $x,y\in K$, let 
$$f(x*y) = \int f(z) d(\delta _x*\delta _y) (z).$$   
For any $\lam \in M^b(K)$ let
$P_\lam $ be the convolution operator on $L^2(K)$ defined by 
$$P_\lam (f) (x) = \int f(x*y) d\lam (y)$$
for all $x \in K$ and for all $f \in L^2(K)$.  It is easy to see that 
for $\mu $ and $\lam$ in $M^b(K)$, $P_{\mu *\lam } = P_\mu 
P_\lam$ and
for $\mu \in {\cal P}(K)$, $P_\mu$ is a contraction and 
$P_{\check \mu} = \check P_\mu $ where $\check P$ is the adjoint of $P$ for any
operator $P$ (see \cite {B-H} for more on convolution operators, for 
$f*\lam$ in 1.2.15 of \cite {B-H}, $P_\lam (f) = f*\check \lam$).  

\br\label{rm1}
Let $M^{(1)}(K)$ be the space all non-negative measures in $M^b(K)$ such 
that $0\leq \mu (K) \leq 1$.  We would like to remark 
following well-known facts regarding convolution operators on $L^2 (K)$.  

\be
\item $\varphi \colon \lam \mapsto P_\lam$ is a Banach algebra 
representation of $M^b(K)$ 
into the space of bounded linear operators on $L^2(K)$.

\item $\varphi (M^{(1)}(K))$ and $\varphi ({\cal P}(K))$ are convex. 

\item $M^{(1)}(K)$ is compact in the vague topology and hence 
$M^{(1)}(K)$ 
with vague topology is isomorphic to $(\varphi (M^{(1)}(K))$ with weak 
operator topology.  In particular, $\varphi (M^{(1)}(K))$ is closed in the 
weak operator topology and hence - being convex - also in the strong 
operator topology.

\item On ${\cal P}(K)$ respectively $\varphi ({\cal P}(K))$ the above 
topologies and weak topology coincide.
\ee
\er

A probability measure $\lam \in {\cal P}(K)$ is called {\it scattered} if
$\sup _{x \in K} \delta _x * \lam ^n (C) \ra 0$ 
for all compact sets $C$ in $K$ where 
$\lam ^n$ is the $n$-th convolution power of $\lam$.  
It is known that for a locally compact group $G$, $\lam \in {\cal P}(G)$
is scattered if and only if $||P_\lam ^n (f)|| \ra 0$ for all $f \in
L^2(G)$ (see Proposition 3.3 of \cite{JRW}).  It may be proved in a
similar way that $\lam$ on a hypergroup is scattered if and only if 
$||P_\lam ^n (f)|| \ra 0$ for all $f \in L^2 (K)$.  

For $\lam , \mu \in M^{(1)}(K)$ and $x \in K$, $\mu\lam$ and $x \mu$ 
denote $\mu *\lam$ and $\delta _x *\mu$.  

Let $E_0 = \{ f\in L^2 (K) \mid ||P_\lam ^n (f) || \ra 0 \}$.  Suppose
that a probability measure $\lam \in {\cal P}(K)$ is not scattered, then
$L^2 (K) \not = E_0$.  
Bartoszek and Rebowski studied all $f \in L^2 (G)$
for which $||P_\lam ^n (f) || \not \ra 0$ for adapted probability measures
$\lam$ (adapted probability measures are those probability measures for 
which the closed subgroup generated by the support is the whole group)
on certain class of groups, namely 
it is proved in \cite{BR} that for
adapted probability measures on groups $G$ with left and right uniform 
structures
are equivalent, $$L^2 (G) = E_0 \oplus L^2 (G, \Sigma _d) \eqno (1)$$
where $\Sigma _d$ is the deterministic $\sigma$-field associated to the
Markov operator $P_\lam$ and if $\lam$ is non-scattered, then $(L^2 (G,
\Sigma _d), P_\lam)$ is isomorphic to the bilateral shift on $l^2 (\mathbb
Z)$.  This type of decomposition of $L^2 (G)$ was first studied in
\cite{KL}
and for compact groups and abelian groups \cite{KL} gives an affirmative
answer.  The decomposition (1) of $L^2(K)$ is known as {\it Krengel-Lin
decomposition}.  Here, we are interested in proving the afore-stated
Krengel-Lin decomposition for probabilities on general hypergroups.

We briefly sketch the results proved in  this article.  In section 2, we
prove that a probability measure $\lam$ on a hypergroup $K$ has
Krengel-Lin
decomposition if and only if $(\check \lam ^n \lam ^n)$ converges to an
idempotent in ${\cal P}(K)$ or $\lam$ is scattered 
(see Theorem \ref{t2,1}) and in the remaining sections we
verify this condition for Tortrat groups (see Theorem \ref{t3,1}), for
commutative hypergroups (see Theorem \ref{t4,1}) and for central
hypergroups (see Theorem \ref{t5,1}).  In section 5, we give an example to
show that Krengel-Lin decomposition does not hold for 
certain measures on discrete hypergroups.
\end{section}

\begin{section}{Markov operators and measures on hypergroups}

Let $(X, \Sigma , m)$ be a $\sigma$-finite measure space. A {\it Markov
operator} on $(X, \Sigma , m)$ is a linear contraction $P\colon L^\infty
(X) \ra L^\infty (X)$ such that $P$ preserves the cone of non-negative
functions, $P(1) =1$ and $P(f_n ) \downarrow 0$ a.e. 
if $f _n \downarrow 0$ and $0\leq f_n \leq 1$ in $L^\infty$.  The measure
$m$ is called {\it invariant} if
$\int P(f) dm = \int f dm$ for all $f$.  In that case $P$ is also a
contraction on $L^1(X)$ and therefore in all spaces $L^p(X)$ for $1\leq
p\leq \infty$ (see \cite{KL}).  The convolution operators of probability 
measures on hypergroups (admitting a invariant measure) are examples of 
Markov operators.

The deterministic $\sigma$-algebra, $\Sigma _d$ associated to $P$ is
defined as the $\sigma$-algebra of measurable sets $A$ in $\Sigma$ such
that for each $n \geq 1$, $P^n (\chi _A ) = \chi _{B_n}$ for some
measurable set $B_n$ in $G$.  The deterministic $\sigma$-algebra was
introduced to study the asymptotic behavior of the iterates of $P$.  

We now recall the following results from \cite{F}: 

\noindent {\bf Theorem F}  Let $P$ be a Markov operator on a
$\sigma$-finite measure space $(X, \Sigma , m)$ with invariant measure
$m$.  Then 
\be
\item [(i)] $L^2 (X, \Sigma _d) = \{ f \in L^2 (X) \mid \check P^n P^n (f)
= f ~~{\rm for ~~ all ~~} n \geq 1 \}$ where $\check P$ is the adjoint of
$P$ on $L^2 (X)$ and 

\item [(ii)] if $f \perp L^2 (X, \Sigma _d)$, then $P^n (f) \ra 0$ in the
weak topology.
\ee

It may be easily seen that the Krengel-Lin decomposition holds for $P$ if
and only if we can have strong convergence in Theorem F(ii).  It is known
that in general we cannot have strong convergence in Theorem F(ii) (see
\cite{KL} and references cited there).  Thus, the Krengel-Lin
decomposition
does not hold for any Markov operator. 

We first state a Proposition for Markov operators on $L^2$-spaces.

\bp\label{pn2,2}
Let $P$ be a Markov operator on a $\sigma$-finite measure space $(X,
\Sigma , m)$ with an invariant measure $m$.  Then  

\be
\item [(i)] there exists a operator $Q$ on $L^2(X)$ such that 
$\check P^n P^n \ra Q$ in the strong operator topology, 

\item [(ii)] if $P$ is a normal operator, then $Q^2 =Q$ and 

\item [(iii)] $\check P^n P^n (f) \ra 0$ weakly 
(hence strongly in view of (i)) implies that $P^n (f) \ra 0$ strongly.  
\ee
\ep

\bo
For a Markov operator $P$, (i) follows from a well-known result 
known as convergence of alternating 
sequences, here we include a proof of it.  The sequence 
$(\check P^n P^n)$ is a decreasing sequence of positive contractions and 
hence $(\check P^n P^n )$ converges in the strong operator topology.  

Suppose $P$ is a normal operator, we have 
$\check P^n P^n = (\check P P)^n$ and hence $\check PP Q=Q$.  
This implies that $Q^2= Q$.  Thus, proving (ii) and that (iii) is easy 
to verify.
\eo

We now deduce from Theorem F and Proposition \ref{pn2,2}, a necessary and 
sufficient condition for a Markov operator to have Krengel-Lin 
decomposition: It may be mentioned that a similar result may be found in 
\cite{Ro} Chapter IV, 4, Lemma 3.

\bc\label{p2,2}
Let $P$ be a Markov operator on a $\sigma$-finite measure space $(X,
\Sigma , m)$ with an invariant measure $m$.  
Then $P$ has Krengel-Lin decomposition if and only
if $\check P^n P^n \ra Q$ in the weak operator topology where $Q$ is the
projection onto $L^2(X, \Sigma _d)$.  
\ec

We now interpret the necessary and sufficient condition in Corollary 
\ref{p2,2} for measures on hypergroups.  

\bt\label{t2,1}
Let $K$ be a hypergroup and $\lam$ be a probability measure on $K$.  Then
$\lam$ has Krengel-Lin decomposition, that is $L^2 (K) = E_0 \oplus L^2
(K, \Sigma _d)$ if and only if either $\lam $ is scattered or $(\check
\lam ^n \lam ^n)$ converges to an idempotent in ${\cal P}(K)$.
\et

\bo
From Corollary \ref{p2,2} and Remark \ref{rm1}, we get 
that $\lam$ has Krengel-Lin decomposition if and only if 
$\check \lam ^n \lam ^n \ra \rho$ in the vague topology where 
$\rho$ is an idempotent in $M^{(1)}(K)$ such that $L^2(K, \Sigma _d)$ is 
the space of all $P_\rho$ fixed functions.  For any idempotent $\rho$ in 
$M^{(1)}(K)$, either $\rho =0$ or $\rho \in {\cal P}(K)$.  
If $\check \lam ^n \lam ^n \ra \rho = \rho ^2 \in {\cal P}(K)$, then 
by Theorem F(i) and since $\check \lam ^n \lam ^n \rho = \rho$, we 
get that $L^2(K, \Sigma _d)$ is the space of all $P_\rho$ fixed 
points.  This shows that $\lam$ has Krengel-Lin decomposition if and only 
if either $\lam$ is scattered or $(\check \lam ^n \lam ^n)$ converges to 
an idempotent.
\eo

\end{section}

\begin{section}{Measures on Tortrat groups}

In this section we prove the Krengel-Lin decomposition for probability
measures on Tortrat groups.  A locally compact group $G$ is called {\it
Tortrat} if for any sequence of the form $(x_n \lam x_n^{-1})$ has an
idempotent limit point in ${\cal P}(G)$ only if $\lam$ is an idempotent.
Tortrat class was introduced by P. Eisele, this class contains all
SIN-groups and all distal linear groups (see \cite{E2} and \cite{R}).  

We now prove the Krengel-Lin decomposition for certain probabilities.  We
first recall that for a locally compact group $G$, a probability measure
$\lam \in {\cal P}(G)$ is called {\it adapted} if the closed subgroup
generated by the support of $\lam$ is $G$ itself.  The structure of
non-scattered adapted probability measures on groups is well-studied in 
\cite{J}.  In \cite {J}, under some additional structural conditions on
$G$ or if $\lam$ is spread-out (a power of $\mu$ is not singular with 
respect to the Haar measure), it is proved that there 
exists a $g \in G$ such that
$(g^{-n}\lam ^n)$ converges but in view of a result in \cite{E1} we are
interested in studying the cases for
which there is a $g \in G$ such that $(g^{-n}\lam ^n)$ converges to an
idempotent.

\bp\label{p3,1}
Let $G$ be a non-compact locally compact group and $\lam$
be an adapted regular Borel probability measure on $G$.  Suppose there
exist a compact normal subgroup $K$ such that $\check \lam ^n \lam ^n \ra
\omega _K$.  Then we have the following:
\be
\item $L^2 (G) = E_0 \oplus L^2 (G, \Sigma _d)$;

\item $\Sigma _d$ is the $\sigma$-algebra generated by $\{x^n K \mid n \in
\mathbb Z\}$ for any $x$ in the support of $\lam$;

\item $(L^2 (G, \Sigma _d), P_\lam )$ is isomorphic to the bilateral shift
on $l^2 (\mathbb Z)$.

\ee
\ep

\bo
Since $\check \lam ^n \lam ^n \ra \omega _K$, (1) follows from Theorem
2.1.

We now claim that $K$ is the smallest closed normal subgroup a coset of
which contains support of $\lam$.  By Theorem 4.3 of \cite{E1}, there
exists a
$x \in G$ such that $x^{-n} \lam ^n \ra \omega _K$.  This implies that
$x^{-1} \omega _K \lam = \omega _K$.  Since $K$ is normal, $\omega _K
x^{-1}\lam = \omega _K$.  This implies that $\lam$ is supported on $Kx$.
It is easy to see that $K$ is contained in any closed normal subgroup a
coset of which contains the support of $\lam$.  This proves the claim.

Now, the rest of proof closely follows \cite{BR}.  Since $\lam$ is
adapted,
$K$ is open and hence by normalizing $m$, we may assume that $\omega _K
(f) = \int fdm$.  Now, for any $x$ in the support of $\lam$, we have
$$P_\lam ^n (\chi _{x^mK} )=\chi _{x^{m-n}K} \eqno (i)$$  for all $m$ and
$n$ in $\mathbb Z$.  This implies that the $\sigma $-algebra generated by
$\{x^n K \mid n \in \mathbb Z \}$ is contained in $\Sigma _d$.  
For $f\in L^2(G, \Sigma _d)$, by Theorem F(i), 
$P_{\check \lam ^n *\lam ^n} (f) =f$ for all 
$n \geq 1$.  Hence by assumption, $P_{\omega _K}(f) =f$.  This implies 
that $f$ is constant on the cosets of $K$.  
Thus, $\Sigma _d$ is the
$\sigma$-algebra generated by $\{x ^n K \mid n \in \mathbb Z \}$ for any
$x$ in the support of $\lam$.  This proves (2) and (3) follows from
equation (i).
\eo

We now prove the Krengel-Lin decomposition for measures on Tortrat groups.

\bt\label{t3,1}
Let $G$ be a non-compact Tortrat group and $\lam$ be an
adapted probability in ${\cal P}(G)$.  Suppose $\lam$ is not scattered.
Then $L^2 (G) = E_0 \oplus L^2 (G, \Sigma _d)$ and $(L^2 (G, \Sigma _d),
P_\lam )$ is isomorphic to the bilateral shift on $l^2 (\mathbb Z)$.
Also, the deterministic $\sigma$-algebra is generated by $\{ g^n K \mid 
n\in \mathbb Z\}$ for any $g$ in the support of $\lam$ and for some
compact normal subgroup $K$ of $G$.
\et

\bo
Since $\lam$ is not scattered, 
$\check \lam ^n \lam ^n \ra \rho \in {\cal P}(G)$.  We first claim that 
$\rho ^2 = \rho$.  Suppose $G$ is metrizable,
by Theorem 1.1 of \cite{E1}, there exists a sequence $(x_n)$ in $G$ such
that
$(x_n \lam ^n)$ converges.  By Theorem 2.1 of \cite{E2}, for all $x$ in
the
support of $\lam$, $x^{-n}\lam ^n \ra \omega _H$ for some compact subgroup
$H$ such that $xH = Hx$.  This implies that $\check \lam ^n \lam ^n \ra
\omega _H$ and $xH=Hx$ for all $x$ in the support of $\lam$ implies 
that $H$ is a normal subgroup since $\lam$ is adapted.  
In the general case, since $\lam$ is 
adapted, $G$ is $\sigma$-compact and hence $G$ can be approximated by 
metrizable groups.  Then by applying standard arguments as in Theorem 3.4 
of \cite {E2}, we prove that $\rho$ is an idempotent.  
Let $K$ be a compact subgroup of $G$ such that $\rho = \omega _K$.  Since 
$\check \lam ^n \rho \lam^n = \rho$ and $\lam$ is adapted we get that 
$K$ is normal in $G$.  Now the result follows from Proposition
\ref{p3,1}.
\eo
\end{section}

\begin{section}{Measures on hypergroups}

In this section we consider probability measures on commutative
hypergroups and central hypergroups.  Krengel-Lin decomposition for 
normal probability measures on hypergroups and hence in particular, for 
probability measures on commutative hypergroups may be easily deduced 
from Proposition \ref{pn2,2} (ii).  

\bt\label{t4,1}
Let $K$ be a hypergroup and $\lam$ be a non-scattered normal probability
measure on $K$.  Then $L^2 (K) = E_0 \oplus L^2 (K, \Sigma _d)$.  In
particular Krengel-Lin decomposition holds for all $\lam \in {\cal P}(K)$
if $K$ is a commutative hypergroup.
\et

We next consider central hypergroups.  
Let $K$ be a hypergroup.  We shall denote the maximal subgroup of $K$ by
$$G(K) = \{ x \in K \mid x*\check x = e \}$$ and the center of $K$ by
$$Z(K) = \{ x \in K \mid x *y = y*x {\rm ~~ for ~~ all~~} y \in K\}.$$
The hypergroup $K$ is called {\it central} if $K/Z$ is compact where $Z =
Z(K) \cap G(K)$; we remark
that $K/Z$ is again a hypergroup.  Central hypergroups arise
naturally as double coset spaces of compact subgroups of central groups
and central hypergroups have invariant Haar measures (see \cite {HKK} for 
a proof of the existence of Haar measures on central hypergroups).    
We first recall the following result on the shift compactness 
of factors and see 5.1.4 of \cite {B-H} for a proof.  

\bp\label{p5,1}
Let $K$ be a metrizable hypergroup.  Let $(\mu _n )$, $(\lam _n)$ and 
$(\eta _n)$ be sequences of probability measures on $K$.  Suppose 
$\mu _n = \eta _n \lam _n $ for all $n \geq 1$ and $(\mu _n)$ is
relatively compact.  Then there exists a 
sequence $(x_n)$ in $K$ such that 
$(x_n \lam _n)$ is relatively compact.
\ep

The following result may be compared with Theorem 3.1 of \cite{C}.

\bp\label{con}
Let $K$ be a metrizable hypergroup and $\lam \in {\cal P}(K)$.  Suppose 
$\lam $ is not scattered.  Then there exists a sequence $(x_n )$ in $K$ 
such that $(x_n \lam ^n )$ is relatively compact.
\ep

\bo
By Proposition \ref{pn2,2}, 
$\check P_\lam ^n P_\lam ^n \ra P$ in the strong operator topology.  
Since $P_\lam ^n (f) \not \ra 0$ for some $f \in L^2 (K)$, 
$P (f) \not = 0$.  By Remark \ref{rm1}, there exists a 
$\rho \in M^{(1)} (K)$ such that $P_\rho = P$ and hence there exists a 
$\rho \in M^{(1)}(K)$ such that 
$\check \lam ^n \rho \lam ^n = \rho$ for all $n \geq 1$.  
Replacing $\rho $ by $\rho /\rho (K)$, we
may assume that there exists a $\rho \in {\cal P} (K)$ such that $\check
\lam ^n \rho \lam ^n = \rho$ for all $n \geq 1$.  By Proposition
\ref{p5,1}, there exists a sequence $(x_n)$ in $K$ such that 
$(x _n \lam ^n)$ is relatively compact.  
\eo

We now prove the Krengel-Lin decomposition for measures on central
hypergroups.  

\bt\label{t5,1}
Let $K$ be a metrizable central hypergroup and $\lam \in {\cal P}(K)$ be
non-scattered.  Then there exists an idempotent $\rho $ such that $\check
\lam ^n \lam ^n \ra \rho $ in ${\cal P}(K)$ and $L^2 (K) = E_0 \oplus
L^2(K, \Sigma _d) $.
\et

\bo
Suppose $\lam$ is not scattered.  By Proposition \ref{con} 
and since $K$ is central hypergroup, there exists
a sequence $(g _n )$ in $Z$ such that $(g_n \lam ^n)$ is relatively
compact.  Then $(\check \lam ^n \lam ^n)$ is relatively compact.  In view 
of Proposition \ref{pn2,2} (i) and Remark \ref{rm1} there exists a 
$\rho \in {\cal P}(K)$ such that $\check \lam ^n \lam ^n\ra \rho$.  

Let $\mu _n = g_{n-1}^{-1} \lam g_n $ and $\nu _k^n = \mu _{k+1} \cdots
\mu _n$ for all $n \geq 1$ and $k <n$ where $g_0 = e$.  Then $(\nu _0^n)$
is relatively compact.  Arguing as in \cite{C}, we may prove that there
exists
a subsequence $n(i)$ such that $\lim \nu _k ^{n(i)} =\tilde \nu _k 
\in {\cal P}(K)$ for all $k \geq 1$ and $\lim \tilde \nu _{n(i)} = \nu
_\infty$.  Also, $\nu _\infty$ is an idempotent in ${\cal P}(K)$.  Now,
for $n >1$, $\lam \lam ^{n-1} g_{n-1}g_{n-1}^{-1} g_n = \lam ^n g_n$.
This implies that $(g_{n-1} ^{-1} g_n)$ is relatively compact.  Now for $n
>1$ and $k <n$, we have 
$\nu _k^n = \lam ^{n-k}g_{n-k}g_{n-k}^{-1}g_k^{-1}g_n$ and 
$\check \nu _k ^n \nu _k^n = \check \lam ^{n-k} \lam ^{n-k}$ 
for all $n >1$ and for all $k < n$.  Thus, $\tilde \nu _k \tilde \nu _k =
\rho$ for all $k \geq 1$.  Thus, $\rho = \nu _\infty$ which is an
idempotent.  This shows that $(\check \lam ^n \lam ^n )$ converges to an
idempotent.  This proves the result.  
\eo

\end{section}

\begin{section}{Example}

It is known that there exists a measure $\lam$ on certain locally
compact groups such that $\lam$ is supported on a coset of a compact 
normal subgroup but $(\check \lam ^n\lam ^n )$ does not converge to
an idempotent and hence $\lam$ does not have Krengel-Lin decomposition
(see \cite{E1} or \cite {J}) .  

We now construct a discrete hypergroup and a measure $\lam$ such 
that $(\check \lam ^n *\lam ^n)$ does not converge to an idempotent.  It 
may remarked that in \cite{Ba} Bartoszek proved that probability measures 
on discrete groups admit Krengel-Lin decomposition but our examples 
show that shifted convolution powers on discrete hypergroups 
as compared to on discrete groups need not have similar behavior.  
Let $A$ be an locally compact abelian group and $\ap$ be an automorphism
of $A$ such that 
\be
\item $\ap (g) \ra e$ for all $g \in A$, 

\item there exists a compact open subgroup $L$ of $A$ and 
$\ap (L) \subset L$ and 

\item there exists a $x \in A$ with $x^6\not \in L$ and 
$\ap (x) \in L$.
\ee

Let $\mu = {1\over 2} (\delta _x +\delta _e) \omega _L$.  
Let $G$ be the semidirect product of $\mathbb Z$ and $A$ where the
$\mathbb Z$ action is given by $\ap$.  
Define $\lam$ by $\check \lam = (1, \mu )$.  
Let $g = (1, e)$.  Then 
$\check \lam ^n g^{-n} = \mu \ap (\mu ) \cdots \ap ^{n-1}(\mu )$.  
It is easy to see that $\mu \ap (\mu ) = \mu$.  Thus, 
$\check \lam ^n g^{-n} = \mu$ for all $n \geq 1$.  Thus, 
$\check \lam ^n \lam ^n \ra \rho$ where 
$\rho = ({1\over 2}\delta _e+{1\over 4}\delta _x +{1\over 4}
\delta _{x^{-1}})\omega _L$.  
Suppose $\rho$ is an idempotent, then either $x^2 $ or $x^3$ is in $L$.  This
implies that $x^6 \in L$ which is a contradiction.  Thus, $\rho$ is not 
an idempotent.  It is easy to see that $\lam$ is $L$-biinvariant.  
Let $K$ be the hypergroup of $L$-double cosets in $G$.  Then $K$ is a
discrete hypergroup and $\lam$  may be viewed as a probability 
measure on $K$.  Thus, we have a discrete hypergroup $K$ and a 
$\lam \in {\cal P}(K)$ such that $(\check \lam ^n \lam ^n )$ does not 
converge to an idempotent.  

We now provide a $A$, $x$, $L$ and $\ap$ satisfying the above
conditions.  Fix a prime integer $p$, let ${\mathbb Q}_p$ be 
the field of $p$-adic numbers and $|\cdot |$ be the $p$-adic
absolute value.  Now take $A = {\mathbb Q}_p$, $x$ with 
$|x| = p^3$, $L= \{ g \in {\mathbb Q}_p \mid |g| \leq 1 \}$ and
$\ap $ is defined by $\ap (g) = p^3 g$ for all $g \in {\mathbb Q}_p$.  
Then $|\ap ^n (g) |=  p^{-3n}|g| \ra 0$ for all $g \in G$.  
Also, $|6x| \geq p^2$ but 
$|\ap (x) |= 1$, that is condition 3 is satisfied and it is easy to check
condition 2.  Using this idea one may construct many such examples.  

\end{section}

\begin{acknowledgement}
I would like to thank the referee of a previous version for fruitful 
suggestions in formulating Proposition \ref{pn2,2} and for many other 
useful comments.
\end{acknowledgement}

\noindent {C. Robinson Edward Raja \\
Stat-Math Unit \\ Indian Statistical Insitute \\
8th Mile Mysore Road \\ R. V. College Post \\
Bangalore 560 059. \\  India. }

\vskip 0.25in

\noindent {creraja@isibang.ac.in}

\end{document}